\documentclass[runningheads]{llncs}

\usepackage{makeidx}  

\usepackage{amsmath}
\usepackage{amssymb}
\usepackage{amstext}
\usepackage{bm}
\usepackage[pdfauthor={Jason Gross, Adam Chlipala, and David I. Spivak},
            pdftitle={Experience Implementing a Performant Category-Theory Library in Coq},
            pdfsubject={Interactive Theorem Proving Conference, 2014},
            pdfkeywords={Coq, category theory, homotopy type theory, duality, performance}]{hyperref}
\usepackage{ifpdf}
\ifpdf
\else
  \usepackage[anythingbreaks]{breakurl}[2013/04/10]
\fi
\usepackage{comment}
\usepackage{xifthen}
\usepackage{longtable}
\usepackage{etex}
\usepackage{xy}
\usepackage{xypic}[2010/04/13]
\usepackage{mathtools}
\usepackage{stmaryrd} 
\usepackage[utf8]{inputenc}
\usepackage{tikz}

\input{glyphtounicode}
\makeatletter
\newlength{\assizeof@height}
\newlength{\assizeof@depth}
\newlength{\assizeof@width}
\newcommand{\asheightof}[2]{%
  \settoheight{\assizeof@height}{#2}%
  \settowidth{\assizeof@width}{#2}%
  \settodepth{\assizeof@depth}{#2}%
  \raisebox{0pt}[\assizeof@height][\assizeof@depth]{#1}%
}
\newsavebox{\checkmarkdashed@box}
\sbox{\checkmarkdashed@box}{%
\begin{tikzpicture}[scale=0.0018, x=\baselineskip, y=-\baselineskip]
\begin{scope}[shift={(-76.14285,-468.51261)}]
  \path[fill=gray] (346.8401,509.3677) .. controls (340.5769,505.0388) and
    (340.7618,500.5122) .. (356.8188,488.4276) .. controls (372.8758,476.3431) and
    (385.5965,469.4556) .. (393.7791,468.6693) .. controls (402.4251,467.8385) and
    (404.1912,473.2846) .. (397.5467,480.2873) .. controls (395.7973,482.1311) and
    (389.1909,487.6189) .. (382.8659,492.4825) .. controls (374.6709,499.4135) and
    (356.1486,515.8013) .. (346.8401,509.3677) -- cycle;
  \path[fill=gray] (260.1089,583.4074) .. controls (266.3140,575.7764) and
    (272.7509,568.2670) .. (279.3743,560.8342) .. controls (324.4398,518.7462) and
    (319.3756,546.8324) .. (289.7561,585.4171) .. controls (282.4707,594.9346) and
    (276.2845,605.1740) .. (270.0874,615.4138) .. controls (215.7257,693.3061) and
    (215.7106,638.2351) .. (260.1089,583.4074) -- cycle;
  \path[fill=gray] (184.4815,709.4074) .. controls (188.7556,700.1304) and
    (193.0479,690.9313) .. (196.6062,682.0703) .. controls (214.1074,660.9549) and
    (227.6958,691.0301) .. (219.6701,711.9379) .. controls (210.8162,735.0030) and
    (203.2310,759.0839) .. (197.0533,784.3044) .. controls (172.5104,824.2648) and
    (170.7368,739.2609) .. (184.4815,709.4074) -- cycle;
  \path[fill=gray] (121.6449,767.9646) .. controls (118.8624,763.3407) and
    (98.3906,735.2020) .. (96.1182,731.3920) .. controls (91.5735,723.7721) and
    (88.0116,717.7146) .. (85.2823,712.9577) .. controls (79.8236,703.4439) and
    (77.6952,699.1324) .. (77.6952,697.9285) .. controls (77.6952,691.5044) and
    (92.7170,679.4441) .. (106.4137,674.8717) .. controls (114.2814,672.2452) and
    (115.9904,672.2686) .. (119.1161,675.0458) .. controls (120.5090,676.2833) and
    (141.7234,710.5494) .. (153.0933,729.2494) .. controls (155.8390,738.5209) and
    (141.8256,804.9784) .. (121.6449,767.9646) -- cycle;
\end{scope}
\end{tikzpicture}%
}

\newcommand{\checkmarkdashed}{\asheightof{\usebox{\checkmarkdashed@box}}{\checkmark}}
\makeatother

\makeatletter
\newcommand{\ensuretext}[1]{%
  \ifmmode
    \expandafter\@firstoftwo
  \else
    \expandafter\@secondoftwo
  \fi
  {\text{#1}}%
  {#1}%
}

\newcommand{\cat}[1]{\ensuremath{\mathcal{#1}}}
\newcommand{\maybesub}[1][]{\ifthenelse{\equal{#1}{}}{}{\ensuremath{_{#1}}}}
\newcommand{\maybecat}[1][]{\ifthenelse{\equal{#1}{}}{}{\ensuremath{_{\cat{#1}}}}}
\newcommand{\Ob}[1][]{\ensuretext{Ob}\maybecat[#1]}
\newcommand{\Hom}[1][]{\ensuretext{Hom}\maybecat[#1]}
\newcommand{\Id}[1][]{\ensuremath{1}\maybesub[#1]}

\DeclareUnicodeCharacter{2218}{\ensuremath{\circ}}
\DeclareUnicodeCharacter{2192}{\ensuremath{\to}}
\DeclareUnicodeCharacter{2200}{\ensuremath{\forall}}
\DeclareUnicodeCharacter{1D52}{\ensuremath{{}^{\text{o}}}}
\DeclareUnicodeCharacter{1D56}{\ensuremath{{}^{\text{p}}}}
\DeclareUnicodeCharacter{03C6}{\ensuremath{\varphi}}
\DeclareUnicodeCharacter{03C8}{\ensuremath{\psi}}
\DeclareUnicodeCharacter{03B1}{\ensuremath{\alpha}}
\DeclareUnicodeCharacter{03B2}{\ensuremath{\beta}}
\DeclareUnicodeCharacter{03B3}{\ensuremath{\gamma}}
\DeclareUnicodeCharacter{03BB}{\ensuremath{\lambda}}
\DeclareUnicodeCharacter{03B7}{\ensuremath{\eta}}
\DeclareUnicodeCharacter{0394}{\ensuremath{\Delta}}
\DeclareUnicodeCharacter{0280}{\ensuremath{{}_0}}
\DeclareUnicodeCharacter{0281}{\ensuremath{{}_1}}
\DeclareUnicodeCharacter{0282}{\ensuremath{{}_2}}
\DeclareUnicodeCharacter{1D30}{\ensuremath{{}^D}}
\DeclareUnicodeCharacter{21D2}{\ensuremath{\Rightarrow}}
\makeatother

\newcommand{\leftfraction}[2]{%
  \makebox[#1\width][l]{#2}%
  \makebox[0pt][l]{%
    \color{white}%
    \rule[-0.5pt]
      {\widthof{#2}}%
      {\totalheightof{#2}+2pt}%
  }%
}

\usepackage{verbatim}
\usepackage{xcolor}

\makeatletter

\begingroup
  \catcode`{=11
  \catcode`}=11
  \catcode`<=1
  \catcode`>=2
  \gdef\@makeopenbrace<\catcode`{>
  \gdef\@makeclosebrace<\catcode`}>
  \global\let\mytextbraceleft={
  \global\let\mytextbraceright=}
\endgroup

\newcommand{\processcommands}[1]{{%
  \catcode`\\=0\relax
  \@makeopenbrace=1\relax
  \@makeclosebrace=2\relax
  \def\{{\texttt{\mytextbraceleft}}
  \def\}{\texttt{\mytextbraceright}}
  \scantokens{#1}%
}}
\newcommand{\processcommandsinverbatimline}{\expandafter\processcommands\expandafter{\the\verbatim@line}}
\newenvironment{coqcode}{\begingroup
  \vspace{0.5\baselineskip}

  \let\item\relax
  \setlength{\parsep}{0pt}%
  \setlength{\parskip}{0pt}%
  \setlength{\topsep}{0pt}%
  \setlength{\@topsepadd}{0pt}%
  \setlength{\partopsep}{0pt}%
  \setlength{\@topsep}{0pt}%
  \let\old@@par\@@par
  \let\@@par\relax
  \let\old@vskip\vskip
  \let\vskip\relax
  \verbatim
  \let\@@par\old@@par
  \let\vskip\old@vskip
  \let\verbatim@processline=\processcommandsinverbatimline
}{\endverbatim\endgroup\vspace{0.5\baselineskip}}
\makeatother

\definecolor{varpurple}{rgb}{0.4,0,0.4}
\definecolor{constrmaroon}{rgb}{0.6,0,0}
\definecolor{defgreen}{rgb}{0,0.4,0}
\definecolor{indblue}{rgb}{0,0,0.8}
\definecolor{kwred}{rgb}{0.8,0.1,0.1}

\newcommand{\colortext}[2]{\textcolor{#1}{#2}}

\newcommand{\coqdockw}[1]{\texttt {\colortext{kwred}{#1}}}

\newcommand{\coqdocvar}[1]{\colortext{varpurple}{#1}}

\newcommand{\coqdoccst}[1]{\texttt{\colortext{defgreen}{#1}}}

\newcommand{\coqdocind}[1]{\texttt{\colortext{indblue}{#1}}}

\newcommand{\coqdocconstr}[1]{\texttt {\colortext{constrmaroon}{#1}}} 

\newcommand{\coqdocinductive}[1]{\coqdocind{#1}}
\newcommand{\coqdocdefinition}[1]{\coqdoccst{#1}}
\newcommand{\coqdocvariable}[1]{\coqdocvar{#1}}
\newcommand{\coqdocconstructor}[1]{\coqdocconstr{#1}}

\newcommand{\coqdocmethod}[1]{\coqdoccst{#1}}

\newcommand{\coqdocrecord}[1]{\coqdocind{#1}}
\newcommand{\coqdocprojection}[1]{\coqdoccst{#1}}
\newcommand{\coqdocnotation}[1]{\coqdockw{#1}}

\newcommand{\pullbacksymbol}{{\ensuremath{\lrcorner}}}
\newcommand{\pullbackarrow}[1][dr]{\ar@{}[#1]|<<\pullbacksymbol}

\usepackage[
  sortcites=true,
  backend=bibtex
]{biblatex}

\DeclareFieldFormat{labelnumberwidth}{#1.}

\begin{document}
\pagestyle{headings}


\mainmatter

\title{Experience Implementing a Performant Category-Theory Library in Coq}

\titlerunning{Experience Implementing a Performant Category-Theory Library in Coq}  
%
\author{Jason Gross \and Adam Chlipala \and David I.~Spivak}
\authorrunning{J.~Gross \and A.~Chlipala \and D.~I.~Spivak} 
%
\tocauthor{Jason Gross, Adam Chlipala, and David I.~Spivak}
\institute{Massachusetts Institute of Technology, Cambridge, MA, USA \\
\email{jgross@mit.edu}, \email{adamc@csail.mit.edu}, \email{dspivak@math.mit.edu}}

\maketitle

\begin{abstract}
  We describe our experience implementing a broad category-theory library in Coq.  Category theory and computational performance are not usually mentioned in the same breath, but we have needed substantial engineering effort to teach Coq to cope with large categorical constructions without slowing proof script processing unacceptably.  In this paper, we share the lessons we have learned about how to represent very abstract mathematical objects and arguments in Coq and how future proof assistants might be designed to better support such reasoning.  One particular encoding trick to which we draw attention allows category-theoretic arguments involving \emph{duality} to be internalized in Coq's logic with definitional equality.  Ours may be the largest Coq development to date that uses the relatively new Coq version developed by homotopy type theorists, and we reflect on which new features were especially helpful.
\end{abstract}




\section{Introduction}

Category theory~\cite{mac1998categories} is a popular all-encompassing mathematical formalism that casts familiar mathematical ideas from many domains in terms of a few unifying concepts.  A \emph{category} can be described as a directed graph plus algebraic laws stating equivalences between paths through the graph.  Because of this spartan philosophical grounding, category theory is sometimes referred to in good humor as ``formal abstract nonsense.''  Certainly the popular perception of category theory is quite far from pragmatic issues of implementation.  This paper is an experience report on an implementation of category theory that has run squarely into issues of design and efficient implementation of type theories, proof assistants, and developments within them.

It would be reasonable to ask, what would it even mean to implement ``formal abstract nonsense,'' and what could the answer have to do with optimized execution engines for functional programming languages?  We mean to cover the whole scope of category theory, which includes many concepts that are not manifestly computational, so it does not suffice merely to employ the well-known folklore semantic connection between categories and typed functional programming~\cite{pierce1988taste}.  Instead, a more appropriate setting is a computer proof assistant.  We chose to build a library for Coq~\cite{coq}, a popular system based on constructive type theory.

One might presume that it is a routine exercise to transliterate categorical concepts from the whiteboard to Coq.  Most category theorists would probably be surprised to learn that standard constructions ``run too slowly,'' but in our experience that is exactly the result of experimenting with na\"ive first Coq implementations of categorical constructs.  It is important to tune the library design to minimize the cost of manipulating terms and proving interesting theorems.

This design experience is also useful for what it reveals about the consequences of design decisions for type theories themselves.  Though type theories are generally simpler than widely used general-purpose programming languages, there is still surprising subtlety behind the few choices that must be made.  Homotopy type theory~\cite{HoTTBook} is a popular subject of study today, where there is intense interest in designing a type theory that makes proofs about topology particularly natural, via altered treatment of equality.  In this setting and others, there remain many open questions about the consequences of type theoretical features for different sorts of formalization.  Category theory, said to be ``notoriously hard to formalize''~\cite{harrison1996formalized}, provides a good stress test of any proof assistant, highlighting problems in usability and efficiency.

Formalizing the connection between universal morphisms and adjunctions provides a typical example of our experience with performance.  A \emph{universal morphism} is a construct in category theory generalizing extrema from calculus.  An \emph{adjunction} is a weakened notion of equivalence.  In the process of rewriting our library to be compatible with homotopy type theory, we discovered that cleaning up this construction conceptually resulted in a significant slow-down, because our first attempted rewrite resulted in a leaky abstraction barrier and, most importantly, large goals (\autoref{sec:term-size}).  Plugging the holes there reduced goal sizes by two orders of magnitude\footnote{The word count of the larger of the two relevant goals went from 163,811 to 1,390.}, which led to a factor of ten speedup in that file (from 39s to 3s), but incurred a factor of three slow-down in the file where we defined the abstraction barriers (from 7s to 21s).  Working around slow projections of $\Sigma$ types (\autoref{sec:records-vs-sigma}) and being more careful about code reuse each gave us back half of that lost time.

For reasons that we present in the course of the paper, we were unsatisfied with the feature set of released Coq version 8.4.  We wound up adopting the Coq version under development by homotopy type theorists~\cite{HoTT/coq}, making critical use of its stronger universe polymorphism (\autoref{sec:category-of-categories}) and higher inductive types (\autoref{sec:hit}).  We hope that our account here provides useful data points for proof assistant designers about which features can have serious impact on proving convenience or performance in very abstract developments.  The two features we mentioned earlier in the paragraph can simplify the Coq user experience dramatically, while a number of other features, at various stages of conception or implementation by Coq team members, can make proving much easier or improve proof script performance by orders of magnitude, generally by reducing term size (\autoref{sec:term-size}): primitive record projections (\autoref{sec:prim-record-proj}), internalized proof irrelevance for equalities (\autoref{sec:equality-reflection}), and $\eta$ rules for records (\autoref{sec:no-judgmental-eta}) and equality proofs (\autoref{sec:compute-match}).

Although pre-existing formalizations of category theory in proof assistants abound~\cite{Ahrens2013,megacz-coq-categories,o2004towards,copumpkin/categories,ConCaT%
,CatsInZFC%
,spitters2010developing%
,kozen2006automating%
,math-overflow-formalizations%
,benediktahrens/coinductives,benediktahrens-Foundations-typesystems,Carvalho1998,jmchapman/restriction-categories,konn/category-agda,crypto-agda/crypto-agda,mattam82-cat,Coalgebras,Algebra,ahrens2010categorical,weber02program,pcapriotti/agda-categories,huet2000constructive,altucher1990mechanically,Category2-AFP,MathClasses,benediktahrens/rezk-completion,mohri1995formalization,spiwackverified%
,caccamo2001higher,aczel1993galois,wilander2005bicategory,logical2001implicit,dyckhoff1985category,wilander2012constructing,harrison1996formalized,agerholm1995experiments,nuo2013second%
}, 
 we chose to implement our library~\cite{HoTT/HoTT-categories} from scratch.  Beginning from scratch allowed the first author to familiarize himself with both category theory and Coq, without simultaneously having to familiarize himself with a large pre-existing code base.  Additionally, starting from scratch forced us to confront all of the decisions involved in designing such a library, and gave us the confidence to change the definitions of basic concepts multiple times to try out various designs, including fully rewriting the library at least three times.  Although this paper is much more about the design of category theory libraries in general than our library in particular, we include a comparison of our library~\cite{HoTT/HoTT-categories} with selected extant category theory libraries in \autoref{sec:compare-libraries}.  At present, our library subsumes many of the constructions in most other such Coq libraries, and is not lacking any constructions in other libraries that are of a complexity requiring significant type checking time, other than monoidal categories.

\medskip

We begin our discussion in \autoref{sec:categories} considering a mundane aspect of type definitions that has large consequences for usability and performance.  With the expressive power of Coq's logic Gallina, we often face a choice of making \emph{parameters} of a type family explicit arguments to it, which looks like universal quantification; or of including them within values of the type, which looks like existential quantification.  As a general principle, we found that the universal or \emph{outside} style improves the user experience modulo performance, while the existential or \emph{inside} style speeds up type checking.  The rule that we settled on was: \emph{inside} definitions for pieces that are usually treated as black boxes by further constructions, and \emph{outside} definitions for pieces whose internal structure is more important later on.

\autoref{sec:duality-unification} presents one of our favorite design patterns for categorical constructions: a way of coaxing Coq's definitional equality into implementing \emph{proof by duality}, one of the most widely known ideas in category theory.  In \autoref{sec:other}, we describe a few other design choices that had large impacts on usability and performance, often of a few orders of magnitude.  \autoref{sec:compare-libraries} wraps up with a grid comparison of our library with others.

\section{Issues in Defining the Type of Categories}\label{sec:categories}
  We have chosen to use the outside style when we care more about the definition of a construct than about carrying it around as an opaque blob to fit into other concepts.  The first example of this choice comes up in deciding how to define categories.

  \subsection{Dependently Typed Morphisms}
    In standard mathematical practice, a category \cat{C} can be defined~\cite{awodey2010category} to consist of:
    \begin{itemize}
      \item
        a class \Ob[C] of \emph{objects}
      \item
        for all objects $a, b \in \Ob[C]$, a class $\Hom[C](a, b)$ of \emph{morphisms from $a$ to $b$}
      \item
        for each object $x \in \Ob[C]$, an \emph{identity morphism} $\Id[x] \in \Hom[C](x, x)$
      \item
        for each triple of objects $a, b, c \in \Ob[C]$, a \emph{composition function} $\circ : \Hom[C](b, c) \times \Hom[C](a, b) \to \Hom[C](a, c)$
    \end{itemize}
    satisfying the following axioms:
    \begin{itemize}
      \item
        associativity: for composable morphisms $f$, $g$, $h$, we have $f \circ (g \circ h) = (f \circ g) \circ h$.
      \item
        identity: for any morphism $f \in \Hom[C](a, b)$, we have $\Id[b] \circ f = f = f \circ \Id[a]$
    \end{itemize}

    Following \cite{HoTTBook}, we additionally require our morphisms to be 0-truncated (to have unique identity proofs).  Without this requirement, we have a standard pre-- homotopy type theory definition of a category.

    We might\footnote{The definition we actually use has some additional fields; see, e.g., \autoref{sec:remove-symmetry}.} formalize the definition in Coq (if Coq had mixfix notation) as: \label{def:category}
\begin{coqcode}
\coqdockw{Record} \coqdocrecord{Category} :=
  \{ \coqdocprojection{Ob} : \coqdockw{Type};
    \coqdocprojection{Hom} : \coqdocmethod{Ob} → \coqdocmethod{Ob} → \coqdockw{Type};
    _\coqdocprojection{∘}_ : \coqdockw{∀} \{\coqdocvariable{a} \coqdocvariable{b} \coqdocvariable{c}\}, \coqdocmethod{Hom} \coqdocvariable{b} \coqdocvariable{c} → \coqdocmethod{Hom} \coqdocvariable{a} \coqdocvariable{b} → \coqdocmethod{Hom} \coqdocvariable{a} \coqdocvariable{c};
    \coqdocprojection{1} : \coqdockw{∀} \{\coqdocvariable{x}\}, \coqdocmethod{Hom} \coqdocvariable{x} \coqdocvariable{x};
    \coqdocprojection{Assoc} : \coqdockw{∀} \coqdocvariable{a} \coqdocvariable{b} \coqdocvariable{c} \coqdocvariable{d} (\coqdocvariable{f} : \coqdocmethod{Hom} \coqdocvariable{c} \coqdocvariable{d}) (\coqdocvariable{g} : \coqdocmethod{Hom} \coqdocvariable{b} \coqdocvariable{c}) (\coqdocvariable{h} : \coqdocmethod{Hom} \coqdocvariable{a} \coqdocvariable{b}),
      \coqdocvariable{f} \coqdocmethod{∘} (\coqdocvariable{g} \coqdocmethod{∘} \coqdocvariable{h}) = (\coqdocvariable{f} \coqdocmethod{∘} \coqdocvariable{g}) \coqdocmethod{∘} \coqdocvariable{h};
    \coqdocprojection{LeftId} : \coqdockw{∀} \coqdocvariable{a} \coqdocvariable{b} (\coqdocvariable{f} : \coqdocmethod{Hom} \coqdocvariable{a} \coqdocvariable{b}), \coqdocmethod{1} \coqdocmethod{∘} \coqdocvariable{f} = \coqdocvariable{f};
    \coqdocprojection{RightId} : \coqdockw{∀} \coqdocvariable{a} \coqdocvariable{b} (\coqdocvariable{f} : \coqdocmethod{Hom} \coqdocvariable{a} \coqdocvariable{b}), \coqdocvariable{f} \coqdocmethod{∘} \coqdocmethod{1} = \coqdocvariable{f};
    \coqdocprojection{Truncated} : \coqdockw{∀} \coqdocvariable{a} \coqdocvariable{b} (\coqdocvariable{f} \coqdocvariable{g} : \coqdocmethod{Hom} \coqdocvariable{a} \coqdocvariable{b}) (\coqdocvariable{p} \coqdocvariable{q} : \coqdocvariable{f} = \coqdocvariable{g}), \coqdocvariable{p} = \coqdocvariable{q} \}.
\end{coqcode}

    We could just as well have replaced the classes $\Hom[C](a, b)$ with a single class of morphisms \Hom[C], together with functions defining the source and target of each morphism.  Then it would be natural to define morphism composition to take a further argument, a proof of equal domain and codomain between the morphisms.  Users of dependent types are aware that explicit manipulation of equality proofs can complicate code substantially, often to the point of obscuring what would be the heart of an argument on paper.  For instance, the algebraic laws associated with categories must be stated with explicit computation of equality proofs, and further constructions only become more involved.  Additionally, such proofs will quickly bloat the types of goals, resulting in slower type checking.  For these reasons, we decided to stick with the definition of \texttt{Category} above, getting more lightweight help from the type checker in place of explicit proofs.

    \subsection{Complications from Categories of Categories}

    Some complications arise in applying the last subsection's definition of categories to the full range of common constructs in category theory.  One particularly prominent example formalizes the structure of a collection of categories, showing that this collection itself may be considered as a category.

    The morphisms in such a category are \emph{functors},\label{sec:define-functor} maps between categories consisting of a function on objects, a function on hom-types, and proofs that these functions respect composition and identity~\cite{mac1998categories,awodey2010category,HoTTBook}.

    The na\"ive concept of a ``category of all categories,'' \label{sec:category-of-categories} which includes even itself, leads into mathematical inconsistencies, which manifest as universe inconsistency errors in Coq.  The standard resolution is to introduce a hierarchy of categories, where, for instance, most intuitive constructions are considered \emph{small} categories, and then we also have \emph{large} categories, one of which is the category of small categories.  Both definitions wind up with literally the same text in Coq, giving:
\begin{coqcode}
\coqdockw{Definition} \coqdocdefinition{SmallCat} : \coqdocrecord{LargeCategory} :=
  \{| \coqdocprojection{Ob} := \coqdocrecord{SmallCategory};
     \coqdocprojection{Hom} \coqdocvariable{C} \coqdocvariable{D} := \coqdocrecord{SmallFunctor} \coqdocvariable{C} \coqdocvariable{D}; ... |\}.
\end{coqcode}

    It seems a shame to copy-and-paste this definition (and those of \texttt{Category}, \texttt{Functor}, etc.) $n$ times to define an $n$-level hierarchy.  Coq 8.4 and some earlier versions support a flavor of \emph{universe polymorphism} that allows the universe of a definition to vary as a function of the universes of its arguments.  Unfortunately, it is not natural to parametrize \texttt{Cat} by anything but a universe level, which does not have first-class status in Coq anyway.  We found the connection between universe polymorphism and arguments to definitions to be rather inconvenient, and it forced us to modify the definition of \texttt{Category} so that the record field \texttt{Ob} changes into a parameter of the type family.  Then we were able to use the following slightly awkward construction:
\begin{coqcode}
\coqdockw{Definition} \coqdocdefinition{Cat_helper} \coqdocvariable{I} \coqdocvariable{ObOf} (\coqdocvariable{CatOf} : \coqdockw{∀} \coqdocvariable{i} : \coqdocvariable{I}, \coqdocrecord{Category} (\coqdocvariable{ObOf} \coqdocvariable{i}))
    : \coqdocrecord{Category} \coqdocvariable{I}
 := \{| \coqdocprojection{Hom} \coqdocvariable{C} \coqdocvariable{D} := \coqdocrecord{Functor} (\coqdocvariable{CatOf} \coqdocvariable{C}) (\coqdocvariable{CatOf} \coqdocvariable{D}); ... |\}.
\coqdockw{Notation} \coqdocdefinition{Cat} := (\coqdocdefinition{Cat_helper} \{\coqdocvariable{T} : \coqdockw{Type} \& \coqdocrecord{Category} \coqdocvariable{T}\} \coqdocprojection{projT1} \coqdocprojection{projT2}).
\end{coqcode}
Now the definition is genuinely reusable for an infinite hierarchy of sorts of category, because the \texttt{Notation} gives us a fresh universe each time we invoke it, but we have paid the price of adding extra parameters to both \texttt{Category} and \texttt{Cat\_helper}, and this seemingly innocent change produces substantial blow-up in the sizes of proof goals arising during interesting constructions.  So, in summary, we decided that the basic type theoretical design of Coq 8.4 did not provide very good support for pleasing definitions that can be reasoned about efficiently.

This realization (and a few more that will come up shortly) pushed us to become early adopters of the modified version of Coq developed by homotopy type theorists~\cite{HoTT/coq}.  Here, an established kind of more general universe polymorphism~\cite{Harper1991107}, previously implemented only in NuPRL, is available, and the definitions we wanted from the start work as desired.

  \subsection{Arguments vs.~Fields} \label{sec:arguments-vs-fields}

    Unlike most of our other choices, there is a range of possibilities in defining categories, with regards to arguments (on the outside) and fields (on the inside).  At one extreme, everything can be made a field, with a type \texttt{Category} whose inhabitants are categories.  At the other extreme, everything can be made an argument to a dummy function.  Some authors \cite{spitters2010developing} have chosen the intermediate option of making all of the computationally relevant parts (objects, morphisms, composition, and the identity morphism) arguments and the irrelevant proofs (associativity and left and right identity) fields.  We discussed in \autoref{sec:category-of-categories} the option of parametrizing on just the type of objects.  We now consider pros and cons of other common options; we found no benefits to the ``outside'' extreme.
    
    \subsubsection{Everything on the Inside}
      Once we moved to using the homotopy type theorists' Coq with its broader universe polymorphism, we decided to use fields for all of the components of a category.  Switching from the version where the types of objects and morphisms were parameters brought a factor of three speed-up in compilation time over our whole development.  The reason is that, at least in Coq, the performance of proof tree manipulations depends critically on their size~(\autoref{sec:term-size}).  By contrast, the size of the normal form of the term does not seem to matter much in most constructions; see \autoref{sec:duality-unification} for an explanation and the one exception that we might have found.  By using fields rather than parameters for the types of objects and morphisms, the type of functors goes from
      \begin{align*}
        \texttt{\coqdocrecord{Functor} : } & \texttt{\coqdockw{$\forall$} (\coqdocvariable{ob$_{\cat C}$} : \coqdockw{Type}) (\coqdocvariable{ob$_{\cat D}$} : \coqdockw{Type})} \\
        & \texttt{(\coqdocvariable{hom$_{\cat C}$} : \coqdocvariable{ob$_{\cat C}$} $\to$ \coqdocvariable{ob$_{\cat C}$} $\to$ \coqdockw{Type}) (\coqdocvariable{hom$_{\cat D}$} : \coqdocvariable{ob$_{\cat D}$} $\to$ \coqdocvariable{ob$_{\cat D}$} $\to$ \coqdockw{Type}),} \\
        & \texttt{\coqdocrecord{Category} \coqdocvariable{ob$_{\cat C}$} \coqdocvariable{hom$_{\cat C}$} $\to$ \coqdocrecord{Category} \coqdocvariable{ob$_{\cat D}$} \coqdocvariable{hom$_{\cat D}$} $\to$ \coqdockw{Type}}
      \end{align*}
      \noindent to
      \begin{center}
        \texttt{\coqdocrecord{Functor} : \coqdocrecord{Category} $\to$ \coqdocrecord{Category} $\to$ \coqdockw{Type}}
      \end{center}
      The corresponding reduction for the type of natural transformations is even more remarkable, and with a construction that uses natural transformations multiple times, the term size blows up very quickly, even with only two parameters.  If we had more parameters (for composition and identity), the term size would blow up even more quickly.

  	  Usually, we do not care what objects and morphisms a particular category has; most of our constructions take as input arbitrary categories.  Thus, there is a significant performance benefit to having all of the fields on the inside and so hidden from most theorem statements.

    \subsubsection{Relevant Things on the Outside}
      One of the main benefits to making all of the relevant components arguments, and requiring all of the fields to satisfy proof irrelevance, is that it allows the use of type-class resolution without having to worry about overlapping instances.  Practically, this choice means that it is easier to get Coq to infer automatically the proofs that given types and operations assemble into a category, at least in simple cases.  Although others~\cite{spitters2010developing} have found this approach useful, we have not found ourselves wishing we had type-class resolution when formalizing constructions, and there is a significant computational cost of exposing so many parameters in types.  The ``packed classes'' of Ssreflect~\cite{Garillot2009%
      } alleviate this problem by combining this approach with the previous one, at the slight cost of more verbosity in initial definitions.

  \subsection{Equality} \label{sec:equality}
    Equality, which has recently become a very hot topic in type theory~\cite{HoTTBook} and higher category theory~\cite{Leinster2007}, provides another example of a design decision where most usage is independent of the exact implementation details.  Although the question of what it means for objects or morphisms to be equal does not come up much in classical 1-category theory, it is more important when formalizing category theory in a proof assistant, for reasons seemingly unrelated to its importance in higher category theory.  We consider some possible notions of equality.

    \subsubsection{Setoids}
      A setoid~\cite{bishop1967foundations} is a carrier type equipped with an equivalence relation; a map of setoids is a function between the carrier types and a proof that the function respects the equivalence relations of its domain and codomain.  Many authors \cite{copumpkin/categories,MathClasses,megacz-coq-categories,huet2000constructive
      }
      choose to use a setoid of morphisms, which allows for the definition of the category of set(oid)s, as well as the category of (small) categories, without assuming functional extensionality, and allows for the definition of categories where the objects are quotient types.  However, there is significant overhead associated with using setoids everywhere, which can lead to slower compile times.  Every type that we talk about needs to come with a relation and a proof that this relation is an equivalence relation.  Every function that we use needs to come with a proof that it sends equivalent elements to equivalent elements.  Even worse, if we need an equivalence relation on the universe of ``types with equivalence relations,'' we need to provide a transport function between equivalent types that respects the equivalence relations of those types.

    \subsubsection{Propositional Equality}
      An alternative to setoids is propositional equality, which carries none of the overhead of setoids, but does not allow an easy formulation of quotient types, and requires assuming functional extensionality to construct the category of sets.

      Intensional type theories like Coq's have a built-in notion of equality, often called definitional equality or judgmental equality, and denoted as $x \equiv y$.  This notion of equality, which is generally internal to an intensional type theory and therefore cannot be explicitly reasoned about inside of that type theory, is the equality that holds between $\beta\delta\iota\zeta\eta$-convertible terms.

      Coq's standard library defines what is called \emph{propositional equality} on top of judgmental equality, denoted $x = y$. One is allowed to conclude that propositional equality holds between any judgmentally equal terms.

      Using propositional equality rather than setoids is convenient because there is already significant machinery made for reasoning about propositional equalities, and there is much less overhead.  However, we ran into significant trouble when attempting to prove that the category of sets has all colimits, which amounts to proving that it is closed under disjoint unions and quotienting; quotient types cannot be encoded without assuming a number of other axioms.

    \subsubsection{Higher Inductive Types}\label{sec:hit}
      The recent emergence of higher inductive types allows the best of both worlds.  The idea of higher inductive types~\cite{HoTTBook} is to allow inductive types to be equipped with extra proofs of equality between constructors. They originated as a way to allow homotopy type theorists to construct types with non-trivial higher paths.  A very simple example is the interval type, from which functional extensionality can be proven~\cite{interval-implies-funext}.\footnote{This assumes a computational interpretation of higher inductives, an open problem.}  The interval type consists of two inhabitants \texttt{zero}~\texttt{:}~\texttt{Interval} and \texttt{one}~\texttt{:}~\texttt{Interval}, and a proof \texttt{seg~:~zero~=~one}.  In a hypothetical type theory with higher inductive types, the type checker does the work of carrying around an equivalence relation on each type for us, and forbids users from constructing functions that do not respect the equivalence relation of any input type.  For example, we can, hypothetically, prove functional extensionality as follows:
\begin{coqcode}
\coqdockw{Definition} \coqdocdefinition{f_equal} \{\coqdocvariable{A B x y}\} (\coqdocvariable{f} : \coqdocvariable{A} → \coqdocvariable{B}) : \coqdocvariable{x} = \coqdocvariable{y} → \coqdocvariable{f x} = \coqdocvariable{f y}.
\coqdockw{Definition} \coqdocdefinition{functional_extensionality} \{\coqdocvariable{A B}\} (\coqdocvariable{f g} : \coqdocvariable{A} → \coqdocvariable{B})
    : (\coqdockw{∀} \coqdocvariable{x}, \coqdocvariable{f x} = \coqdocvariable{g x}) → \coqdocvariable{f} = \coqdocvariable{g}
  := \coqdockw{λ} (\coqdocvariable{H} : \coqdockw{∀} \coqdocvariable{x}, \coqdocvariable{f x} = \coqdocvariable{g x})
       ⇒ \coqdocdefinition{f_equal} (\coqdockw{λ} (\coqdocvariable{i} : \coqdocind{Interval}) (\coqdocvariable{x} : \coqdocvariable{A})
                     ⇒ \coqdockw{match} \coqdocvariable{i} \coqdockw{with}
                         | \coqdocconstructor{zero} ⇒ \coqdocvariable{f x}
                         | \coqdocconstructor{one}  ⇒ \coqdocvariable{g x}
                         | \coqdocconstructor{seg}  ⇒ \coqdocvariable{H x}
                       \coqdockw{end})
                  \coqdocconstructor{seg}.
\end{coqcode}
      Had we neglected to include the branch for \texttt{seg}, the type checker should complain about an incomplete match; the function \texttt{\coqdockw{$\lambda$}~\coqdocvariable{i}~:~\coqdocinductive{Interval} $\Rightarrow$~\coqdockw{match}~\coqdocvariable{i}~\coqdockw{with} \coqdocconstructor{zero}~$\Rightarrow$~\coqdocconstructor{true} |~\coqdocconstructor{one}~$\Rightarrow$~\coqdocconstructor{false} \coqdockw{end}} of type \texttt{Interval $\to$ bool} should not typecheck for this reason.

      The key insight is that most types do not need any special equivalence relation, and, moreover, if we are not explicitly dealing with a type with a special equivalence relation, then it is impossible (by parametricity) to fail to respect the equivalence relation.  Said another way, the only way to construct a function that might fail to respect the equivalence relation would be by some eliminator like pattern matching, so all we have to do is guarantee that direct invocations of the eliminator result in functions that respect the equivalence relation.

      As with the choice involved in defining categories, using propositional equality with higher inductive types rather than setoids derives many of its benefits from not having to deal with all of the overhead of custom equivalence relations in constructions that do not need them.  In this case, we avoid the overhead by making the type checker or the metatheory deal with the parts we usually do not care about.  Most of our definitions do not need custom equivalence relations, so the overhead of using setoids would be very large for very little gain.  We plan to use higher inductive types\footnote{We fake these in Coq using Yves Bertot's Private Inductive Types extension~\cite{Bertot2013}.} to define quotients, which are necessary to show the existence of certain functors involving the category of sets.  We also currently use higher inductive types to define propositional truncation~\cite{HoTTBook}, which we use to define what it means for a function to be surjective, and prove that in the category of sets, being an isomorphism (an invertible morphism) is equivalent to being injective and surjective.



  \subsection{Records vs.~Nested \texorpdfstring{$\mathrm{\Sigma}$}{Σ} Types} \label{sec:records-vs-sigma}
    In Coq, there are two ways to represent a data structure with one constructor and many fields: as a single inductive type with one constructor (records), or as a nested $\Sigma$ type.  For instance, consider a record type with two type fields $A$ and $B$ and a function $f$ from $A$ to $B$.  A logically equivalent encoding would be $\Sigma A. \; \Sigma B. \; A \to B$.  There are two important differences between these encodings in Coq.

    \label{sec:prim-record-proj}
    The first is that while a theorem statement may abstract over all possible $\Sigma$ types, it may not abstract over all record types, which somehow have a less first-class status.  Such a limitation is inconvenient and leads to code duplication.

    The far more pressing problem, overriding the previous point, is that nested $\Sigma$ types have horrendous performance, and are sometimes a few orders of magnitude slower.  The culprit is projections from nested $\Sigma$ types, which, when unfolded (as they must be, to do computation), each take almost the entirety of the nested $\Sigma$ type as an argument, and so grow in size very quickly.  Matthieu Sozeau is currently working on primitive projections for records for Coq, which would eliminate this problem by eliminating the arguments to the projection functions.\footnote{We eagerly await the day when we can take advantage of this feature in our library.}

\section{Internalizing Duality Arguments in Type Theory} \label{sec:duality-unification}
  In general, we have tried to design our library so that trivial proofs on paper remain trivial when formalized.  One of Coq's main tools to make proofs trivial is the definitional equality, where some facts follow by computational reduction of terms.  We came up with some small tweaks to core definitions that allow a common family of proofs by \emph{duality} to follow by computation.

  Proof by duality is a common idea in higher mathematics: sometimes, it is productive to flip the directions of all the arrows.  For example, if some fact about least upper bounds is provable, chances are that the same kind of fact about greatest lower bounds will also be provable in roughly the same way, by replacing ``greater than''s with ``less than''s and vice versa.

  Concretely, there is a dualizing operation on categories that inverts the directions of the morphisms:
\begin{coqcode}
\coqdockw{Notation} "C ᵒᵖ" := (\{| \coqdocprojection{Ob} := \coqdocprojection{Ob} \coqdocvariable{C}; \coqdocprojection{Hom} \coqdocvariable{x y} := \coqdocprojection{Hom} \coqdocvariable{C y x}; ... |\}).
\end{coqcode}

  Dualization can be used, roughly, for example, to turn a proof that Cartesian product is an associative operation into a proof that disjoint union is an associative operation; products are dual to disjoint unions.

  One of the simplest examples of duality in category theory is initial and terminal objects.  In a category \cat C, an initial object $0$ is one that has a unique morphism $0 \to x$ to every object $x$ in \cat C; a terminal object $1$ is one that has a unique morphism $x \to 1$ from every object $x$ in \cat C.  Initial objects in \cat C are terminal objects in $\cat{C}^\text{op}$.  The initial object of any category is unique up to isomorphism; for any two initial objects $0$ and $0'$, there is an isomorphism $0 \cong 0'$.  By flipping all of the arrows around, we can prove, by duality, that the terminal object is unique up to isomorphism.  More precisely, from a proof that an initial object of $\cat{C}^{\text{op}}$ is unique up to isomorphism, we get that any two terminal objects $1'$ and $1$ in $\cat{C}$, which are initial in $\cat{C}^{\text{op}}$, are isomorphic in $\cat{C}^{\text{op}}$.  Since an isomorphism $x \cong y$ in $\cat{C}^\text{op}$ is an isomorphism $y \cong x$ in \cat C, we get that $1$ and $1'$ are isomorphic in $\cat C$.

  It is generally straightforward to see that there is an isomorphism between a theorem and its dual, and the technique of dualization is well-known to category theorists, among others.  We discovered that, by being careful about how we defined things, we could make theorems be judgmentally equal to their duals!  That is, when we prove a theorem
  \begin{align*}
  \texttt{\coqdocdefinition{initial\_ob\_unique} : \coqdockw{$\forall$} } & \texttt{\coqdocvariable{C} (\coqdocvariable{x y} : \coqdocprojection{Ob} \coqdocvariable{C}),} \\ & \texttt{\coqdocdefinition{is\_initial\_ob} \coqdocvariable{x} $\to$ \coqdocdefinition{is\_initial\_ob} \coqdocvariable{y} $\to$ \coqdocvariable{x} $\cong$ \coqdocvariable{y}},
  \end{align*}
  we can define another theorem
  \begin{align*}
  \texttt{\coqdocdefinition{terminal\_ob\_unique} : \coqdockw{$\forall$} } & \texttt{\coqdocvariable{C} (\coqdocvariable{x y} : \coqdocprojection{Ob} \coqdocvariable{C}),} \\ & \texttt{\coqdocdefinition{is\_terminal\_ob} \coqdocvariable{x} $\to$ \coqdocdefinition{is\_terminal\_ob} \coqdocvariable{y} $\to$ \coqdocvariable{x} $\cong$ \coqdocvariable{y}}
  \end{align*}
  as
  \begin{center}
  \texttt{\coqdocdefinition{terminal\_ob\_unique} \coqdocvariable{C x y H H'} := \coqdocdefinition{initial\_ob\_unique} \coqdocvariable{C}\coqdocnotation{$^{\texttt{op}}$} \coqdocvariable{y x H' H}}.
  \end{center}
  Interestingly, we found that in proofs with sufficiently complicated types, it can take a few seconds or more for Coq to accept such a definition; we are not sure whether this is due to peculiarities of the reduction strategy of our version of Coq, or speed dependency on the size of the normal form of the type (rather than on the size of the unnormalized type), or something else entirely.

  In contrast to the simplicity of witnessing the isomorphism, it takes a significant amount of care in defining concepts, often to get around deficiencies of Coq, to achieve \emph{judgmental} duality.  Even now, we were unable to achieve this ideal for some theorems.  For example, category theorists typically identify the functor category $\cat{C}^\text{op} \to \cat{D}^\text{op}$ (whose objects are functors $\cat{C}^\text{op} \to \cat{D}^\text{op}$ and whose morphisms are natural transformations) with $(\cat{C} \to \cat{D})^\text{op}$ (whose objects are functors $\cat{C} \to \cat{D}$ and whose morphisms are flipped natural transformations).  These categories are canonically isomorphic (by the dualizing natural transformations), and, with the univalence axiom~\cite{HoTTBook}, they are equal as categories!  But we have not found a way to make them definitionally equal, much to our disappointment.

  \subsection{Duality Design Patterns}
    One of the simplest theorems about duality is that it is involutive; we have that $(\cat{C}^{\text{op}})^{\text{op}} = \cat{C}$.  In order to internalize proof by duality via judgmental equality, we sometimes need this equality to be judgmental.  Although it is impossible in general in Coq 8.4 (see \hyperref[sec:no-judgmental-eta]{dodging judgmental \texorpdfstring{$\eta$}{η} on records} below), we want at least to have it be true for any explicit category (that is, any category specified by giving its objects, morphisms, etc., rather than referred to via a local variable).

    \subsubsection{Removing Symmetry} \label{sec:remove-symmetry}
      Taking the dual of a category, one constructs a proof that $f \circ (g \circ h) = (f \circ g) \circ h$ from a proof that $(f \circ g) \circ h = f \circ (g \circ h)$.  The standard approach is to apply symmetry.  However, because applying symmetry twice results in a judgmentally different proof, we decided instead to extend the definition of \texttt{Category} to require both a proof of $f \circ (g \circ h) = (f \circ g) \circ h$ and a proof of $(f \circ g) \circ h = f \circ (g \circ h)$.  Then our dualizing operation simply swaps the proofs.  We added a convenience constructor for categories that asks only for one of the proofs, and applies symmetry to get the other one.  Because we formalized 0-truncated category theory, where the type of morphisms is required to have unique identity proofs, asking for this other proof does not result in any coherence issues.

    \subsubsection{Dualizing the Terminal Category}
      To make everything work out nicely, we needed the terminal category, which is the category with one object and only the identity morphism, to be the dual of itself.  We originally had the terminal category as a special case of the discrete category on $n$ objects.  Given a type $T$ with uniqueness of identity proofs, the discrete category on $T$ has as objects inhabitants of $T$, and has as morphisms from $x$ to $y$ proofs that $x = y$.  These categories are not judgmentally equal to their duals, because the type $x = y$ is not judgmentally the same as the type $y = x$.  As a result, we instead used the indiscrete category, which has \texttt{unit} as its type of morphisms.

    \subsubsection{Which Side Does the Identity Go On?}
      The last tricky obstacle we encountered was that when defining a functor out of the terminal category, it is necessary to pick whether to use the right identity law or the left identity law to prove that the functor preserves composition; both will prove that the identity composed with itself is the identity.  The problem is that dualizing the functor leads to a road block where either concrete choice turns out to be ``wrong,'' because the dual of the functor out of the terminal category will not be judgmentally equal to another instance of itself.  To fix this problem, we further extended the definition of category to require a proof that the identity composed with itself is the identity.

    \subsubsection{Dodging Judgmental \texorpdfstring{$\eta$}{η} on Records}  \label{sec:no-judgmental-eta}
      The last problem we ran into was the fact that sometimes, we really, really wanted judgmental $\eta$ on records.  The $\eta$ rule for records says any application of the record constructor to all the projections of an object yields exactly that object; e.g. for pairs, $x \equiv (x_1, x_2)$ (where $x_1$ and $x_2$ are the first and second projections, respectively).  For categories, the $\eta$ rule says that given a category \cat C, for a ``new'' category defined by saying that its objects are the objects of \cat C, its morphisms are the morphisms of \cat C, \ldots, the ``new'' category is judgmentally equal to \cat C.

      In particular, we wanted to show that any functor out of the terminal category is the opposite of some other functor; namely, any $F : 1 \to \cat C$ should be equal to $(F^{\text{op}})^{\text{op}} : 1 \to (\cat C^{\text{op}})^{\text{op}}$.  However, without the judgmental $\eta$ rule for records, a local variable $\cat C$ cannot be judgmentally equal to $(\cat C^{\text{op}})^{\text{op}}$, which reduces to an application of the constructor for a category.  To get around the problem, we made two variants of dual functors: given $F : \cat C \to \cat D$, we have $F^{\text{op}} : \cat C^{\text{op}} \to \cat D^{\text{op}}$, and given $F : C^{\text{op}} \to \cat D^{\text{op}}$, we have $F^{\text{op}'} : \cat C \to \cat D$.  There are two other flavors of dual functors, corresponding to the other two pairings of ${}^{\text{op}}$ with domain and codomain, but we have been glad to avoid defining them so far.  As it was, we ended up having four variants of dual natural transformation, and are very glad that we did not need sixteen.  We look forward to Coq 8.5, when we will hopefully only need one.

  \subsection{Moving Forward: Computation Rules for Pattern Matching} \label{sec:compute-match}
    While we were able to work around most of the issues that we had in internalizing proof by duality, things would have been far nicer if we had more $\eta$ rules.  The $\eta$ rule for records is explained above.  The $\eta$ rule for equality says that the identity function is judgmentally equal to the function $f : \forall x\, y, x = y \to x = y$ defined by pattern matching on the first proof of equality; this rule is necessary to have any hope that applying symmetry twice is judgmentally the identity transformation.  Matthieu Sozeau is currently working on giving Coq judgmental $\eta$ for records with one or more fields, though not for equality~\cite{mattam82/coq-polyproj}%
    .

    \autoref{sec:associators} will give more examples of the pain of manipulating pattern matching on equality.  Homotopy type theory provides a framework that systematizes reasoning about proofs of equality, turning a seemingly impossible task into a manageable one.  However, there is still a significant burden associated with reasoning about equalities, because so few of the rules are judgmental.

    We are currently attempting to divine the appropriate computation rules for pattern matching constructs, in the hopes of making reasoning with proofs of equality more pleasant.\footnote{See \url{https://coq.inria.fr/bugs/show\_bug.cgi?id=3179} and \url{https://coq.inria.fr/bugs/show\_bug.cgi?id=3119}.}

\section{Other Design Choices}\label{sec:other}

A few other pervasive strategies made non-trivial differences for proof performance or simplicity.

  \subsection{Identities vs.~Equalities; Associators} \label{sec:associators}
    There are a number of constructions that are provably equal, but which we found more convenient to construct transformations between instead, despite the increased verbosity of such definitions.  This is especially true of constructions that strayed towards higher category theory.  For example, when constructing the Grothendieck construction of a functor to the category of categories, we found it easier to first generalize the construction from functors to pseudofunctors.  The definition of a pseudofunctor results from replacing various equalities in the definition of a functor with isomorphisms (analogous to bijections between sets or types), together with proofs that the isomorphisms obey various coherence properties.  This replacement helped because there are fewer operations on isomorphisms (namely, just composition and inverting), and more operations on proofs of equality (pattern matching, or anything definable via induction); when we were forced to perform all of the operations in the same way, syntactically, it was easier to pick out the operations and reason about them.

    Another example was defining the (co)unit of adjunction composition, where instead of a proof that $F \circ (G \circ H) = (F \circ G) \circ H$, we used a natural transformation, a coherent mapping between the actions of functors.  Where equality-based constructions led to computational reduction getting stuck at casts, the constructions with natural transformations reduce in all of the expected contexts.

  \subsection{Opacity; Linear Dependence of Speed on Term Size}\label{sec:equality-reflection}\label{sec:term-size}

    Coq is slow at dealing with large terms.  For goals around 175,000 words long\footnote{When we had objects as arguments rather than fields (see \autoref{sec:arguments-vs-fields}), we encountered goals of about 219,633 words when constructing pointwise Kan extensions.}, we have found that simple tactics like \texttt{apply f\_equal} take around 1 second to execute, which makes interactive theorem proving very frustrating.\footnote{See also \url{https://coq.inria.fr/bugs/show\_bug.cgi?id=3280}.}  Even more frustrating is the fact that the largest contribution to this size is often arguments to irrelevant functions, i.e., functions that are provably equal to all other functions of the same type.  (These are proofs related to algebraic laws like associativity, carried inside many constructions.)

    Opacification helps by preventing the type checker from unfolding some definitions, but it is not enough: the type checker still has to deal with all of the large arguments to the opaque function.  Hash-consing might fix the problem completely.

    Alternatively, it would be nice if, given a proof that all of the inhabitants of a type were equal, we could forget about terms of that type, so that their sizes would not impose any penalties on term manipulation.  
    One solution might be irrelevant fields, like those of Agda, or implemented via the Implicit CiC~\cite{barras2008implicit,logical2001implicit}.

  \subsection{Abstraction Barriers} \label{sec:abstraction-barriers}

    In many projects, choosing the right abstraction barriers is essential to reducing mistakes, improving maintainability and readability of code, and cutting down on time wasted by programmers trying to hold too many things in their heads at once.  This project was no exception; we developed an allergic reaction to constructions with more than four or so arguments, after making one too many mistakes in defining limits and colimits.  Limits are a generalization, to arbitrary categories, of subsets of Cartesian products.  Colimits are a generalization, to arbitrary categories, of disjoint unions modulo equivalence relations.

    Our original flattened definition of limits involved a single definition with 14 nested binders for types and algebraic properties.  After a particularly frustrating experience hunting down a mistake in one of these components, we decided to factor the definition into a larger number of simpler definitions, including familiar categorical constructs like terminal objects and comma categories.  This refactoring paid off even further when some months later we discovered the universal morphism definition of adjoint functors~\cite{wiki:adjoint-functors:universal-morphisms%
    ,ncatlab:adjoint+functor:UniversalArrows%
    }%
    .  With a little more abstraction, we were able to reuse the same decomposition to prove the equivalence between universal morphisms and adjoint functors, with minimal effort.

    Perhaps less typical of programming experience, we found that picking the right abstraction barriers could drastically reduce compile time by keeping details out of sight in large goal formulas.  In the instance discussed in the introduction, we got a factor of ten speed-up by plugging holes in a leaky abstraction barrier!\footnote{See \url{https://github.com/HoTT/HoTT/commit/eb0099005171} for the exact change.}

\section{Comparison of Category-Theory Libraries} \label{sec:compare-libraries}

\newcommand{\allstars}{**********}

We present here a table comparing the features of various category-theory libraries.  Our library is the first column.  Gray dashed check-marks (\checkmarkdashed) indicate features in progress.  The right-most library is in Agda; the rest are in Coq.  A check-mark with $n$ stars (*) indicates a construction taking 20$n$ seconds to compile on a 64-bit server with a 2.40 GHz CPU and 16 GB of RAM.

\begin{center}
\begin{longtable}{l|ccccc}
\multicolumn{1}{c|}{Construction} & \textbf{\cite{HoTT/HoTT-categories}} & \cite{megacz-coq-categories} & \cite{ConCaT} & \cite{Ahrens2013} & \cite{copumpkin/categories} \\
\hline
\endhead
Mostly automated (with custom Ltac) & \phantom{\leftfraction{0.15}{\allstars}}\checkmark\makebox[0pt][l]{\relax}\phantom{\leftfraction{0.15}{\allstars}} & \phantom{\leftfraction{0.15}{\allstars}}\makebox[0pt][l]{\relax}\phantom{\leftfraction{0.15}{\allstars}} & \phantom{\leftfraction{0.15}{\allstars}}\makebox[0pt][l]{\relax}\phantom{\leftfraction{0.15}{\allstars}} & \phantom{\leftfraction{0.15}{\allstars}}\makebox[0pt][l]{\relax}\phantom{\leftfraction{0.15}{\allstars}} & \phantom{\leftfraction{0.15}{\allstars}}\makebox[0pt][l]{\relax}\phantom{\leftfraction{0.15}{\allstars}} \\
Uses HoTT & \phantom{\leftfraction{0.15}{\allstars}}\checkmark\makebox[0pt][l]{\relax}\phantom{\leftfraction{0.15}{\allstars}} & \phantom{\leftfraction{0.15}{\allstars}}\makebox[0pt][l]{\relax}\phantom{\leftfraction{0.15}{\allstars}} & \phantom{\leftfraction{0.15}{\allstars}}\makebox[0pt][l]{\relax}\phantom{\leftfraction{0.15}{\allstars}} & \phantom{\leftfraction{0.15}{\allstars}}\checkmark\makebox[0pt][l]{\relax}\phantom{\leftfraction{0.15}{\allstars}} & \phantom{\leftfraction{0.15}{\allstars}}\makebox[0pt][l]{\relax}\phantom{\leftfraction{0.15}{\allstars}} \\
Uses type classes & \phantom{\leftfraction{0.15}{\allstars}}\makebox[0pt][l]{\relax}\phantom{\leftfraction{0.15}{\allstars}} & \phantom{\leftfraction{0.15}{\allstars}}\checkmark\makebox[0pt][l]{\relax}\phantom{\leftfraction{0.15}{\allstars}} & \phantom{\leftfraction{0.15}{\allstars}}\makebox[0pt][l]{\relax}\phantom{\leftfraction{0.15}{\allstars}} & \phantom{\leftfraction{0.15}{\allstars}}\makebox[0pt][l]{\relax}\phantom{\leftfraction{0.15}{\allstars}} & \phantom{\leftfraction{0.15}{\allstars}}\makebox[0pt][l]{\relax}\phantom{\leftfraction{0.15}{\allstars}} \\
Setoid of morphisms & \phantom{\leftfraction{0.15}{\allstars}}\makebox[0pt][l]{\relax}\phantom{\leftfraction{0.15}{\allstars}} & \phantom{\leftfraction{0.15}{\allstars}}\checkmark\makebox[0pt][l]{\relax}\phantom{\leftfraction{0.15}{\allstars}} & \phantom{\leftfraction{0.15}{\allstars}}\checkmark\makebox[0pt][l]{\relax}\phantom{\leftfraction{0.15}{\allstars}} & \phantom{\leftfraction{0.15}{\allstars}}\makebox[0pt][l]{\relax}\phantom{\leftfraction{0.15}{\allstars}} & \phantom{\leftfraction{0.15}{\allstars}}\checkmark\makebox[0pt][l]{\relax}\phantom{\leftfraction{0.15}{\allstars}} \\
Uses higher inductive types & \phantom{\leftfraction{0.15}{\allstars}}\checkmark\makebox[0pt][l]{\relax}\phantom{\leftfraction{0.15}{\allstars}} & \phantom{\leftfraction{0.15}{\allstars}}\makebox[0pt][l]{\relax}\phantom{\leftfraction{0.15}{\allstars}} & \phantom{\leftfraction{0.15}{\allstars}}\makebox[0pt][l]{\relax}\phantom{\leftfraction{0.15}{\allstars}} & \phantom{\leftfraction{0.15}{\allstars}}\makebox[0pt][l]{\relax}\phantom{\leftfraction{0.15}{\allstars}} & \phantom{\leftfraction{0.15}{\allstars}}\makebox[0pt][l]{\relax}\phantom{\leftfraction{0.15}{\allstars}} \\
Assumes UIP or equivalent & \phantom{\leftfraction{0.15}{\allstars}}\makebox[0pt][l]{\relax}\phantom{\leftfraction{0.15}{\allstars}} & \phantom{\leftfraction{0.15}{\allstars}}\makebox[0pt][l]{\relax}\phantom{\leftfraction{0.15}{\allstars}} & \phantom{\leftfraction{0.15}{\allstars}}\makebox[0pt][l]{\relax}\phantom{\leftfraction{0.15}{\allstars}} & \phantom{\leftfraction{0.15}{\allstars}}\makebox[0pt][l]{\relax}\phantom{\leftfraction{0.15}{\allstars}} & \phantom{\leftfraction{0.15}{\allstars}}\checkmark\makebox[0pt][l]{\relax}\phantom{\leftfraction{0.15}{\allstars}} \\
Category of sets & \phantom{\leftfraction{0.15}{\allstars}}\checkmark\makebox[0pt][l]{\relax}\phantom{\leftfraction{0.15}{\allstars}} & \phantom{\leftfraction{0.15}{\allstars}}\makebox[0pt][l]{\relax}\phantom{\leftfraction{0.15}{\allstars}} & \phantom{\leftfraction{0.15}{\allstars}}\checkmark\makebox[0pt][l]{\relax}\phantom{\leftfraction{0.15}{\allstars}} & \phantom{\leftfraction{0.15}{\allstars}}\checkmark\makebox[0pt][l]{\relax}\phantom{\leftfraction{0.15}{\allstars}} & \phantom{\leftfraction{0.15}{\allstars}}\checkmark\makebox[0pt][l]{\relax}\phantom{\leftfraction{0.15}{\allstars}} \\
Initial/Terminal objects & \phantom{\leftfraction{0.15}{\allstars}}\checkmark\makebox[0pt][l]{\relax}\phantom{\leftfraction{0.15}{\allstars}} & \phantom{\leftfraction{0.15}{\allstars}}\checkmark\makebox[0pt][l]{\relax}\phantom{\leftfraction{0.15}{\allstars}} & \phantom{\leftfraction{0.15}{\allstars}}\checkmark\makebox[0pt][l]{\relax}\phantom{\leftfraction{0.15}{\allstars}} & \phantom{\leftfraction{0.15}{\allstars}}\checkmark\makebox[0pt][l]{\relax}\phantom{\leftfraction{0.15}{\allstars}} & \phantom{\leftfraction{0.15}{\allstars}}\checkmark\makebox[0pt][l]{\relax}\phantom{\leftfraction{0.15}{\allstars}} \\
(co)limits & \phantom{\leftfraction{0.15}{\allstars}}\checkmark\makebox[0pt][l]{\relax}\phantom{\leftfraction{0.15}{\allstars}} & \phantom{\leftfraction{0.15}{\allstars}}\makebox[0pt][l]{\relax}\phantom{\leftfraction{0.15}{\allstars}} & \phantom{\leftfraction{0.15}{\allstars}}\checkmark\makebox[0pt][l]{\relax}\phantom{\leftfraction{0.15}{\allstars}} & \phantom{\leftfraction{0.15}{\allstars}}\checkmark\makebox[0pt][l]{\relax}\phantom{\leftfraction{0.15}{\allstars}} & \phantom{\leftfraction{0.15}{\allstars}}\checkmark\makebox[0pt][l]{\relax}\phantom{\leftfraction{0.15}{\allstars}} \\
(co)limit functor & \phantom{\leftfraction{0.15}{\allstars}}\checkmark\makebox[0pt][l]{\relax}\phantom{\leftfraction{0.15}{\allstars}} & \phantom{\leftfraction{0.15}{\allstars}}\makebox[0pt][l]{\relax}\phantom{\leftfraction{0.15}{\allstars}} & \phantom{\leftfraction{0.15}{\allstars}}\makebox[0pt][l]{\relax}\phantom{\leftfraction{0.15}{\allstars}} & \phantom{\leftfraction{0.15}{\allstars}}\makebox[0pt][l]{\relax}\phantom{\leftfraction{0.15}{\allstars}} & \phantom{\leftfraction{0.15}{\allstars}}\makebox[0pt][l]{\relax}\phantom{\leftfraction{0.15}{\allstars}} \\
(co)limit adjoint to $\Delta$ & \phantom{\leftfraction{0.15}{\allstars}}\checkmark\makebox[0pt][l]{\relax}\phantom{\leftfraction{0.15}{\allstars}} & \phantom{\leftfraction{0.15}{\allstars}}\makebox[0pt][l]{\relax}\phantom{\leftfraction{0.15}{\allstars}} & \phantom{\leftfraction{0.15}{\allstars}}\makebox[0pt][l]{\relax}\phantom{\leftfraction{0.15}{\allstars}} & \phantom{\leftfraction{0.15}{\allstars}}\makebox[0pt][l]{\relax}\phantom{\leftfraction{0.15}{\allstars}} & \phantom{\leftfraction{0.15}{\allstars}}\makebox[0pt][l]{\relax}\phantom{\leftfraction{0.15}{\allstars}} \\
Fully faithful functors & \phantom{\leftfraction{0.15}{\allstars}}\checkmark\makebox[0pt][l]{\relax}\phantom{\leftfraction{0.15}{\allstars}} & \phantom{\leftfraction{0.15}{\allstars}}\makebox[0pt][l]{\relax}\phantom{\leftfraction{0.15}{\allstars}} & \phantom{\leftfraction{0.15}{\allstars}}\makebox[0pt][l]{\relax}\phantom{\leftfraction{0.15}{\allstars}} & \phantom{\leftfraction{0.15}{\allstars}}\checkmark\makebox[0pt][l]{\relax}\phantom{\leftfraction{0.15}{\allstars}} & \phantom{\leftfraction{0.15}{\allstars}}\checkmark\makebox[0pt][l]{\relax}\phantom{\leftfraction{0.15}{\allstars}} \\
Essentially surjective functors & \phantom{\leftfraction{0.15}{\allstars}}\checkmark\makebox[0pt][l]{\relax}\phantom{\leftfraction{0.15}{\allstars}} & \phantom{\leftfraction{0.15}{\allstars}}\checkmark\makebox[0pt][l]{\relax}\phantom{\leftfraction{0.15}{\allstars}} & \phantom{\leftfraction{0.15}{\allstars}}\makebox[0pt][l]{\relax}\phantom{\leftfraction{0.15}{\allstars}} & \phantom{\leftfraction{0.15}{\allstars}}\checkmark\makebox[0pt][l]{\relax}\phantom{\leftfraction{0.15}{\allstars}} & \phantom{\leftfraction{0.15}{\allstars}}\checkmark\makebox[0pt][l]{\relax}\phantom{\leftfraction{0.15}{\allstars}} \\
Unit-Counit Adjunctions & \phantom{\leftfraction{0.15}{\allstars}}\checkmark\makebox[0pt][l]{\relax}\phantom{\leftfraction{0.15}{\allstars}} & \phantom{\leftfraction{0.15}{\allstars}}\checkmark\makebox[0pt][l]{\relax}\phantom{\leftfraction{0.15}{\allstars}} & \phantom{\leftfraction{0.15}{\allstars}}\checkmark\makebox[0pt][l]{\relax}\phantom{\leftfraction{0.15}{\allstars}} & \phantom{\leftfraction{0.15}{\allstars}}\checkmark\makebox[0pt][l]{\relax}\phantom{\leftfraction{0.15}{\allstars}} & \phantom{\leftfraction{0.15}{\allstars}}\checkmark\makebox[0pt][l]{\relax}\phantom{\leftfraction{0.15}{\allstars}} \\
Hom Adjunctions & \phantom{\leftfraction{0.15}{\allstars}}\checkmark\makebox[0pt][l]{\relax}\phantom{\leftfraction{0.15}{\allstars}} & \phantom{\leftfraction{0.15}{\allstars}}\makebox[0pt][l]{\relax}\phantom{\leftfraction{0.15}{\allstars}} & \phantom{\leftfraction{0.15}{\allstars}}\checkmark\makebox[0pt][l]{\relax}\phantom{\leftfraction{0.15}{\allstars}} & \phantom{\leftfraction{0.15}{\allstars}}\makebox[0pt][l]{\relax}\phantom{\leftfraction{0.15}{\allstars}} & \phantom{\leftfraction{0.15}{\allstars}}\makebox[0pt][l]{\relax}\phantom{\leftfraction{0.15}{\allstars}} \\
Universal morphism adjunctions & \phantom{\leftfraction{0.15}{\allstars}}\checkmark\makebox[0pt][l]{\relax}\phantom{\leftfraction{0.15}{\allstars}} & \phantom{\leftfraction{0.15}{\allstars}}\makebox[0pt][l]{\relax}\phantom{\leftfraction{0.15}{\allstars}} & \phantom{\leftfraction{0.15}{\allstars}}\checkmark\makebox[0pt][l]{\relax}\phantom{\leftfraction{0.15}{\allstars}} & \phantom{\leftfraction{0.15}{\allstars}}\makebox[0pt][l]{\relax}\phantom{\leftfraction{0.15}{\allstars}} & \phantom{\leftfraction{0.15}{\allstars}}\makebox[0pt][l]{\relax}\phantom{\leftfraction{0.15}{\allstars}} \\
Adjoint composition laws & \phantom{\leftfraction{0.15}{\allstars}}\checkmark\makebox[0pt][l]{\leftfraction{0.2172}{\allstars}}\phantom{\leftfraction{0.15}{\allstars}} & \phantom{\leftfraction{0.15}{\allstars}}\makebox[0pt][l]{\relax}\phantom{\leftfraction{0.15}{\allstars}} & \phantom{\leftfraction{0.15}{\allstars}}\makebox[0pt][l]{\relax}\phantom{\leftfraction{0.15}{\allstars}} & \phantom{\leftfraction{0.15}{\allstars}}\makebox[0pt][l]{\relax}\phantom{\leftfraction{0.15}{\allstars}} & \phantom{\leftfraction{0.15}{\allstars}}\checkmark\makebox[0pt][l]{\leftfraction{0.06829}{\allstars} \footnotemark}\footnotetext{The use of proof-irrelevant fields speeds up this construction significantly in Agda.}\phantom{\leftfraction{0.15}{\allstars}} \\
Monoidal categories & \phantom{\leftfraction{0.15}{\allstars}}\checkmarkdashed\makebox[0pt][l]{\relax}\phantom{\leftfraction{0.15}{\allstars}} & \phantom{\leftfraction{0.15}{\allstars}}\checkmark\makebox[0pt][l]{\leftfraction{2.0335}{\allstars}}\phantom{\leftfraction{0.15}{\allstars}} & \phantom{\leftfraction{0.15}{\allstars}}\makebox[0pt][l]{\relax}\phantom{\leftfraction{0.15}{\allstars}} & \phantom{\leftfraction{0.15}{\allstars}}\makebox[0pt][l]{\relax}\phantom{\leftfraction{0.15}{\allstars}} & \phantom{\leftfraction{0.15}{\allstars}}\checkmark\makebox[0pt][l]{\leftfraction{0.686205}{\allstars}}\phantom{\leftfraction{0.15}{\allstars}} \\
Enriched categories & \phantom{\leftfraction{0.15}{\allstars}}\checkmarkdashed\makebox[0pt][l]{\relax}\phantom{\leftfraction{0.15}{\allstars}} & \phantom{\leftfraction{0.15}{\allstars}}\checkmark\makebox[0pt][l]{\leftfraction{0.336}{\allstars}}\phantom{\leftfraction{0.15}{\allstars}} & \phantom{\leftfraction{0.15}{\allstars}}\makebox[0pt][l]{\relax}\phantom{\leftfraction{0.15}{\allstars}} & \phantom{\leftfraction{0.15}{\allstars}}\makebox[0pt][l]{\relax}\phantom{\leftfraction{0.15}{\allstars}} & \phantom{\leftfraction{0.15}{\allstars}}\checkmarkdashed\makebox[0pt][l]{\relax}\phantom{\leftfraction{0.15}{\allstars}} \\
2-categories & \phantom{\leftfraction{0.15}{\allstars}}\makebox[0pt][l]{\relax}\phantom{\leftfraction{0.15}{\allstars}} & \phantom{\leftfraction{0.15}{\allstars}}\makebox[0pt][l]{\relax}\phantom{\leftfraction{0.15}{\allstars}} & \phantom{\leftfraction{0.15}{\allstars}}\makebox[0pt][l]{\relax}\phantom{\leftfraction{0.15}{\allstars}} & \phantom{\leftfraction{0.15}{\allstars}}\makebox[0pt][l]{\relax}\phantom{\leftfraction{0.15}{\allstars}} & \phantom{\leftfraction{0.15}{\allstars}}\checkmark\makebox[0pt][l]{\leftfraction{0.06736}{\allstars}}\phantom{\leftfraction{0.15}{\allstars}} \\
Category of (strict) categories & \phantom{\leftfraction{0.15}{\allstars}}\checkmark\makebox[0pt][l]{\relax}\phantom{\leftfraction{0.15}{\allstars}} & \phantom{\leftfraction{0.15}{\allstars}}\makebox[0pt][l]{\relax}\phantom{\leftfraction{0.15}{\allstars}} & \phantom{\leftfraction{0.15}{\allstars}}\checkmark\makebox[0pt][l]{\relax}\phantom{\leftfraction{0.15}{\allstars}} & \phantom{\leftfraction{0.15}{\allstars}}\makebox[0pt][l]{\relax}\phantom{\leftfraction{0.15}{\allstars}} & \phantom{\leftfraction{0.15}{\allstars}}\checkmark\makebox[0pt][l]{\relax}\phantom{\leftfraction{0.15}{\allstars}} \\
Hom functor & \phantom{\leftfraction{0.15}{\allstars}}\checkmark\makebox[0pt][l]{\relax}\phantom{\leftfraction{0.15}{\allstars}} & \phantom{\leftfraction{0.15}{\allstars}}\checkmark\makebox[0pt][l]{\relax}\phantom{\leftfraction{0.15}{\allstars}} & \phantom{\leftfraction{0.15}{\allstars}}\checkmark\makebox[0pt][l]{\relax}\phantom{\leftfraction{0.15}{\allstars}} & \phantom{\leftfraction{0.15}{\allstars}}\checkmark\makebox[0pt][l]{\relax}\phantom{\leftfraction{0.15}{\allstars}} & \phantom{\leftfraction{0.15}{\allstars}}\checkmark\makebox[0pt][l]{\relax}\phantom{\leftfraction{0.15}{\allstars}} \\
Profunctors & \phantom{\leftfraction{0.15}{\allstars}}\checkmark\makebox[0pt][l]{\relax}\phantom{\leftfraction{0.15}{\allstars}} & \phantom{\leftfraction{0.15}{\allstars}}\makebox[0pt][l]{\relax}\phantom{\leftfraction{0.15}{\allstars}} & \phantom{\leftfraction{0.15}{\allstars}}\makebox[0pt][l]{\relax}\phantom{\leftfraction{0.15}{\allstars}} & \phantom{\leftfraction{0.15}{\allstars}}\makebox[0pt][l]{\relax}\phantom{\leftfraction{0.15}{\allstars}} & \phantom{\leftfraction{0.15}{\allstars}}\checkmark\makebox[0pt][l]{\relax}\phantom{\leftfraction{0.15}{\allstars}} \\
Pseudofunctors & \phantom{\leftfraction{0.15}{\allstars}}\checkmark\makebox[0pt][l]{\leftfraction{0.12315}{\allstars}}\phantom{\leftfraction{0.15}{\allstars}} & \phantom{\leftfraction{0.15}{\allstars}}\makebox[0pt][l]{\relax}\phantom{\leftfraction{0.15}{\allstars}} & \phantom{\leftfraction{0.15}{\allstars}}\makebox[0pt][l]{\relax}\phantom{\leftfraction{0.15}{\allstars}} & \phantom{\leftfraction{0.15}{\allstars}}\makebox[0pt][l]{\relax}\phantom{\leftfraction{0.15}{\allstars}} & \phantom{\leftfraction{0.15}{\allstars}}\makebox[0pt][l]{\relax}\phantom{\leftfraction{0.15}{\allstars}} \\
Kan extensions & \phantom{\leftfraction{0.15}{\allstars}}\checkmark\makebox[0pt][l]{\relax}\phantom{\leftfraction{0.15}{\allstars}} & \phantom{\leftfraction{0.15}{\allstars}}\makebox[0pt][l]{\relax}\phantom{\leftfraction{0.15}{\allstars}} & \phantom{\leftfraction{0.15}{\allstars}}\makebox[0pt][l]{\relax}\phantom{\leftfraction{0.15}{\allstars}} & \phantom{\leftfraction{0.15}{\allstars}}\checkmark\makebox[0pt][l]{\relax}\phantom{\leftfraction{0.15}{\allstars}} & \phantom{\leftfraction{0.15}{\allstars}}\checkmark\makebox[0pt][l]{\relax}\phantom{\leftfraction{0.15}{\allstars}} \\
Pointwise Kan extensions & \phantom{\leftfraction{0.15}{\allstars}}\checkmarkdashed\makebox[0pt][l]{\relax}\phantom{\leftfraction{0.15}{\allstars}} & \phantom{\leftfraction{0.15}{\allstars}}\makebox[0pt][l]{\relax}\phantom{\leftfraction{0.15}{\allstars}} & \phantom{\leftfraction{0.15}{\allstars}}\makebox[0pt][l]{\relax}\phantom{\leftfraction{0.15}{\allstars}} & \phantom{\leftfraction{0.15}{\allstars}}\makebox[0pt][l]{\relax}\phantom{\leftfraction{0.15}{\allstars}} & \phantom{\leftfraction{0.15}{\allstars}}\makebox[0pt][l]{\relax}\phantom{\leftfraction{0.15}{\allstars}} \\
${\cat C^{\cat D}}^{\cat E} \cong \cat C^{\cat D \times \cat E}$; $(\cat C \times \cat D)^{\cat E} \cong \cat C^{\cat E}\times \cat D^{\cat E}$ & \phantom{\leftfraction{0.15}{\allstars}}\checkmark\makebox[0pt][l]{\leftfraction{0.17125}{\allstars}}\phantom{\leftfraction{0.15}{\allstars}} & \phantom{\leftfraction{0.15}{\allstars}}\makebox[0pt][l]{\relax}\phantom{\leftfraction{0.15}{\allstars}} & \phantom{\leftfraction{0.15}{\allstars}}\makebox[0pt][l]{\relax}\phantom{\leftfraction{0.15}{\allstars}} & \phantom{\leftfraction{0.15}{\allstars}}\makebox[0pt][l]{\relax}\phantom{\leftfraction{0.15}{\allstars}} & \phantom{\leftfraction{0.15}{\allstars}}\makebox[0pt][l]{\relax}\phantom{\leftfraction{0.15}{\allstars}} \\
Adjoint Functor Theorem & \phantom{\leftfraction{0.15}{\allstars}}\makebox[0pt][l]{\relax}\phantom{\leftfraction{0.15}{\allstars}} & \phantom{\leftfraction{0.15}{\allstars}}\makebox[0pt][l]{\relax}\phantom{\leftfraction{0.15}{\allstars}} & \phantom{\leftfraction{0.15}{\allstars}}\checkmark\makebox[0pt][l]{\leftfraction{0.01635}{\allstars}}\phantom{\leftfraction{0.15}{\allstars}} & \phantom{\leftfraction{0.15}{\allstars}}\makebox[0pt][l]{\relax}\phantom{\leftfraction{0.15}{\allstars}} & \phantom{\leftfraction{0.15}{\allstars}}\makebox[0pt][l]{\relax}\phantom{\leftfraction{0.15}{\allstars}} \\
Yoneda & \phantom{\leftfraction{0.15}{\allstars}}\checkmark\makebox[0pt][l]{\relax}\phantom{\leftfraction{0.15}{\allstars}} & \phantom{\leftfraction{0.15}{\allstars}}\makebox[0pt][l]{\relax}\phantom{\leftfraction{0.15}{\allstars}} & \phantom{\leftfraction{0.15}{\allstars}}\checkmark\makebox[0pt][l]{\relax}\phantom{\leftfraction{0.15}{\allstars}} & \phantom{\leftfraction{0.15}{\allstars}}\checkmark\makebox[0pt][l]{\relax}\phantom{\leftfraction{0.15}{\allstars}} & \phantom{\leftfraction{0.15}{\allstars}}\checkmark\makebox[0pt][l]{\leftfraction{0.323305}{\allstars}}\phantom{\leftfraction{0.15}{\allstars}} \\
dep. product (oplax lim $F : \cat C \to \text{Cat}$) & \phantom{\leftfraction{0.15}{\allstars}}\checkmark\makebox[0pt][l]{\relax}\phantom{\leftfraction{0.15}{\allstars}} & \phantom{\leftfraction{0.15}{\allstars}}\makebox[0pt][l]{\relax}\phantom{\leftfraction{0.15}{\allstars}} & \phantom{\leftfraction{0.15}{\allstars}}\makebox[0pt][l]{\relax}\phantom{\leftfraction{0.15}{\allstars}} & \phantom{\leftfraction{0.15}{\allstars}}\makebox[0pt][l]{\relax}\phantom{\leftfraction{0.15}{\allstars}} & \phantom{\leftfraction{0.15}{\allstars}}\makebox[0pt][l]{\relax}\phantom{\leftfraction{0.15}{\allstars}} \\
dep. sum (oplax colim $F : \cat C \to \text{Cat}$) & \phantom{\leftfraction{0.15}{\allstars}}\checkmark\makebox[0pt][l]{\leftfraction{0.1308}{\allstars}}\phantom{\leftfraction{0.15}{\allstars}} & \phantom{\leftfraction{0.15}{\allstars}}\makebox[0pt][l]{\relax}\phantom{\leftfraction{0.15}{\allstars}} & \phantom{\leftfraction{0.15}{\allstars}}\makebox[0pt][l]{\relax}\phantom{\leftfraction{0.15}{\allstars}} & \phantom{\leftfraction{0.15}{\allstars}}\makebox[0pt][l]{\relax}\phantom{\leftfraction{0.15}{\allstars}} & \phantom{\leftfraction{0.15}{\allstars}}\checkmark\makebox[0pt][l]{\leftfraction{0.10708}{\allstars}}\phantom{\leftfraction{0.15}{\allstars}} \\
\texttt{(\_/\_)} functor $(\cat C^{\cat A})^{\text{op}}\times \cat C^{\cat B} \to \text{Cat}_{/ \cat A \times \cat B}$ & \phantom{\leftfraction{0.15}{\allstars}}\checkmark\makebox[0pt][l]{\leftfraction{0.6074}{\allstars}}\phantom{\leftfraction{0.15}{\allstars}} & \phantom{\leftfraction{0.15}{\allstars}}\makebox[0pt][l]{\relax}\phantom{\leftfraction{0.15}{\allstars}} & \phantom{\leftfraction{0.15}{\allstars}}\makebox[0pt][l]{\relax}\phantom{\leftfraction{0.15}{\allstars}} & \phantom{\leftfraction{0.15}{\allstars}}\makebox[0pt][l]{\relax}\phantom{\leftfraction{0.15}{\allstars}} & \phantom{\leftfraction{0.15}{\allstars}}\makebox[0pt][l]{\relax}\phantom{\leftfraction{0.15}{\allstars}} \\
Rezk completion & \phantom{\leftfraction{0.15}{\allstars}}\makebox[0pt][l]{\relax}\phantom{\leftfraction{0.15}{\allstars}} & \phantom{\leftfraction{0.15}{\allstars}}\makebox[0pt][l]{\relax}\phantom{\leftfraction{0.15}{\allstars}} & \phantom{\leftfraction{0.15}{\allstars}}\makebox[0pt][l]{\relax}\phantom{\leftfraction{0.15}{\allstars}} & \phantom{\leftfraction{0.15}{\allstars}}\checkmark\makebox[0pt][l]{\relax}\phantom{\leftfraction{0.15}{\allstars}} & \phantom{\leftfraction{0.15}{\allstars}}\makebox[0pt][l]{\relax}\phantom{\leftfraction{0.15}{\allstars}} \\
\hline Mean lines per file & \phantom{\leftfraction{0.15}{\allstars}}78\makebox[0pt][l]{\relax}\phantom{\leftfraction{0.15}{\allstars}} & \phantom{\leftfraction{0.15}{\allstars}}126\makebox[0pt][l]{\relax}\phantom{\leftfraction{0.15}{\allstars}} & \phantom{\leftfraction{0.15}{\allstars}}133\makebox[0pt][l]{\relax}\phantom{\leftfraction{0.15}{\allstars}} & \phantom{\leftfraction{0.15}{\allstars}}407\makebox[0pt][l]{\relax}\phantom{\leftfraction{0.15}{\allstars}} & \phantom{\leftfraction{0.15}{\allstars}}98\makebox[0pt][l]{\relax}\phantom{\leftfraction{0.15}{\allstars}} \\
Total compilation time & \phantom{\leftfraction{0.15}{\allstars}}490s\makebox[0pt][l]{\relax}\phantom{\leftfraction{0.15}{\allstars}} & \phantom{\leftfraction{0.15}{\allstars}}517s\makebox[0pt][l]{\relax}\phantom{\leftfraction{0.15}{\allstars}} & \phantom{\leftfraction{0.15}{\allstars}}21s\makebox[0pt][l]{\relax}\phantom{\leftfraction{0.15}{\allstars}} & \phantom{\leftfraction{0.15}{\allstars}}62s\makebox[0pt][l]{\relax \footnotemark}\footnotetext{Nearly 75\% of the time in this library is spent on properties of functor composition.  Nearly 50\% of this time is spent closing sections, for an as-yet unknown reason.}\phantom{\leftfraction{0.15}{\allstars}} & \phantom{\leftfraction{0.15}{\allstars}}717s\makebox[0pt][l]{\relax}\phantom{\leftfraction{0.15}{\allstars}} \\
Total time w/o monoidal & \phantom{\leftfraction{0.15}{\allstars}}490s\makebox[0pt][l]{\relax}\phantom{\leftfraction{0.15}{\allstars}} & \phantom{\leftfraction{0.15}{\allstars}}43s\makebox[0pt][l]{\relax}\phantom{\leftfraction{0.15}{\allstars}} & \phantom{\leftfraction{0.15}{\allstars}}21s\makebox[0pt][l]{\relax}\phantom{\leftfraction{0.15}{\allstars}} & \phantom{\leftfraction{0.15}{\allstars}}62s\makebox[0pt][l]{\relax}\phantom{\leftfraction{0.15}{\allstars}} & \phantom{\leftfraction{0.15}{\allstars}}579s\makebox[0pt][l]{\relax}\phantom{\leftfraction{0.15}{\allstars}} \\
Median file compilation time & \phantom{\leftfraction{0.15}{\allstars}}0.3s\makebox[0pt][l]{\relax}\phantom{\leftfraction{0.15}{\allstars}} & \phantom{\leftfraction{0.15}{\allstars}}0.4s\makebox[0pt][l]{\relax}\phantom{\leftfraction{0.15}{\allstars}} & \phantom{\leftfraction{0.15}{\allstars}}0.1s\makebox[0pt][l]{\relax}\phantom{\leftfraction{0.15}{\allstars}} & \phantom{\leftfraction{0.15}{\allstars}}0.9s\makebox[0pt][l]{\relax}\phantom{\leftfraction{0.15}{\allstars}} & \phantom{\leftfraction{0.15}{\allstars}}1.5s\makebox[0pt][l]{\relax}\phantom{\leftfraction{0.15}{\allstars}} \\
Total number of files & \phantom{\leftfraction{0.15}{\allstars}}147\makebox[0pt][l]{\relax}\phantom{\leftfraction{0.15}{\allstars}} & \phantom{\leftfraction{0.15}{\allstars}}36\makebox[0pt][l]{\relax}\phantom{\leftfraction{0.15}{\allstars}} & \phantom{\leftfraction{0.15}{\allstars}}105\makebox[0pt][l]{\relax}\phantom{\leftfraction{0.15}{\allstars}} & \phantom{\leftfraction{0.15}{\allstars}}13\makebox[0pt][l]{\relax}\phantom{\leftfraction{0.15}{\allstars}} & \phantom{\leftfraction{0.15}{\allstars}}143\makebox[0pt][l]{\relax}\phantom{\leftfraction{0.15}{\allstars}} \\
Total number of definitions & \phantom{\leftfraction{0.15}{\allstars}}578\makebox[0pt][l]{\relax}\phantom{\leftfraction{0.15}{\allstars}} & \phantom{\leftfraction{0.15}{\allstars}}214\makebox[0pt][l]{\relax}\phantom{\leftfraction{0.15}{\allstars}} & \phantom{\leftfraction{0.15}{\allstars}}995\makebox[0pt][l]{\relax}\phantom{\leftfraction{0.15}{\allstars}} & \phantom{\leftfraction{0.15}{\allstars}}367\makebox[0pt][l]{\relax}\phantom{\leftfraction{0.15}{\allstars}} & \phantom{\leftfraction{0.15}{\allstars}}396\makebox[0pt][l]{\relax}\phantom{\leftfraction{0.15}{\allstars}}
\end{longtable}
\end{center}

\vspace*{\dimexpr-2.5\baselineskip\relax}

In summary, our library includes many constructions from past formalizations, plus a few rather complex new ones.  We test the limits of Coq by applying mostly automated Ltac proofs for these constructions, taking advantage of ideas from homotopy type theory and extensions built to support such constructions. In most cases, we found that term size had the biggest impact on speed. We have summarized our observations on using new features from that extension and on other hypothetical features that could make an especially big difference in our development, and we hope these observations can help guide the conversation on the design of future versions of Coq and other proof assistants.

\paragraph{Acknowledgments.} This work was supported in part by the MIT bigdata@CSAIL initiative, NSF grant CCF-1253229, ONR grant N000141310260, and AFOSR grant FA9550-14-1-0031.  We also thank Benedikt Ahrens, Daniel R.~Grayson, Robert Harper, Bas Spitters, and Edward Z.~Yang for feedback on this paper.


\nocite{*}

\clearpage

\printbibliography

\end{document}